\documentclass{amsart}

\usepackage[top=1in, bottom=1in, left=1in, right=1in]{geometry}
\usepackage{lineno}
\usepackage[english]{babel}
\usepackage{xcolor}
\usepackage{framed}
\usepackage[normalem]{ulem}
\usepackage{mathtools}
\usepackage{amsmath}
\usepackage{amsthm}
\usepackage{amssymb}
\usepackage{enumerate}
\usepackage{stmaryrd}
\usepackage{hyperref}
\usepackage{url}
\usepackage{latexsym}
\usepackage{indentfirst}
\usepackage{placeins}
\usepackage{booktabs}
\usepackage{array, makecell}
\usepackage{colortbl}
\usepackage[capitalize]{cleveref}

\theoremstyle{plain}
\newtheorem{theorem}{Theorem}[section]
\newtheorem{lemma}[theorem]{Lemma}

\newtheorem{cor}[theorem]{Corollary}

\theoremstyle{definition}
\newtheorem{defn}[theorem]{Definition}

\theoremstyle{remark}

\crefname{defn}{Definition}{Definitions}
\crefname{cor}{Corollary}{Corollaries}
\crefname{prop}{Proposition}{Propositions}
\crefname{rmk}{Remark}{Remarks}
\crefname{conj}{Conjecture}{Conjectures}

\newcommand{\abs}[1]{\left|#1\right|}

\begin{document}

\title[Matrix displacement ranks over commutative rings]{Matrix Displacement Ranks and Levinson-Type Recursions over Commutative Rings}


\author{Amelia Craciun}
\address{Department of Mathematics, University of California, Berkeley, CA}
\curraddr{}
\email{acraciun@berkeley.edu}
\thanks{}

\author{Rahul Sarkar}
\address{Department of Mathematics, University of California, Berkeley, CA 94720}
\curraddr{}
\email{rsarkar@berkeley.edu}
\thanks{}



\date{}

\dedicatory{}

\commby{}

\begin{abstract}
Displacement rank measures how far a matrix is from being Toeplitz. We develop this theory over commutative rings. Over any commutative ring, the least number of products of triangular Toeplitz matrices needed to represent a matrix as a sum equals the inner rank of its displacement. Over elementary divisor rings with Smith normal form unique up to units, this inner rank is the number of nonzero invariant factors, and inversion exchanges the positive and negative displacement ranks; hence inverses of invertible Toeplitz matrices again have displacement rank at most two. For Hankel and cyclic displacement we give explicit reconstruction formulas, including a solvability criterion and the general solution in the cyclic case, and we diagonalize circulant matrices by the discrete Fourier transform when $n$ is a unit and the ring contains a principal $n$-th root of unity. Finally, we prove Levinson-type recursions for symmetric block Toeplitz matrices and nonsymmetric scalar Toeplitz matrices, assuming that the leading principal submatrices (at block boundaries, in the block case) are invertible. The resulting solver for nonsymmetric systems uses $\mathcal{O}(n^2)$ ring operations.
\end{abstract}

\maketitle

\section{Introduction}
\label{sec:intro}

Let $Z$ denote the $n\times n$ nilpotent lower-shift matrix. For a square
matrix $A$, the \emph{displacements}
\[
\nabla_{+}(A)=A-ZAZ^{\top}
\qquad\text{and}\qquad
\nabla_{-}(A)=A-Z^{\top}AZ
\]
compare $A$ with a copy of itself shifted one step along the diagonal. A matrix is Toeplitz
exactly when $\nabla_{+}(A)$ vanishes outside its first row and first column, so the rank of
$\nabla_{\pm}(A)$ measures how far $A$ is from being Toeplitz. Kailath, Kung, and Morf
\cite{KKM79} showed that, over the real or complex numbers, this rank is also the smallest
number of terms in a representation of $A$ as a sum of products of a lower- and an
upper-triangular Toeplitz matrix, and that inversion exchanges the two displacement ranks:
$\alpha_{+}(A)=\alpha_{-}(A^{-1})$. In particular, the inverse of an invertible Toeplitz
matrix is a sum of at most two such products, even though it is generally not Toeplitz.
These facts underlie a large family of fast algorithms for structured linear systems.

The classical theory is developed over a field, and its proofs use field properties at
several points. The rank of a matrix is unambiguous; the identity
$\operatorname{rank}(I-BC)=\operatorname{rank}(I-CB)$ is available for square matrices $B,C$ of the same size; and every nonzero
pivot can be divided by. None of these is automatic over a commutative ring $R$. There are
several inequivalent notions of rank over $R$ \cite{Brown93}, and when $R$ has zero divisors
the number of nonzero invariant factors of a matrix cannot be read off from its determinantal
ideals. A nonzero element need
not be a unit, so a recursion that is well defined over the fraction field may fail to be
defined over $R$ itself. Yet structured matrices with entries in rings such as $\mathbb{Z}$,
$\mathbb{Z}/m\mathbb{Z}$, or polynomial rings arise naturally in exact computation, and
Bareiss already considered Toeplitz systems over integral domains \cite{Bareiss72}. This
paper determines which parts of the displacement theory and of the Levinson recursions
remain valid over a commutative ring, and under which hypotheses on the ring and on the
matrix.

\subsection{Main results}
\label{ssec:main-results}

Our results fall into four groups. We describe them informally here. All quantities and
notions used below are defined precisely later in the paper, at the point where they are first
needed, and each result is followed by a reference to its formal statement.

\smallskip
\noindent\textit{Displacement rank as inner rank.}
We measure displacement factorizations by the \emph{inner rank} $\rho$, the least $k$ for
which a matrix factors through $R^{k}$ (\cref{def:inner-rank}). Since the displacement
ranks $\alpha_{\pm}$ are themselves defined as minimal lengths of sums, inner rank is the
natural notion. With it, the Kailath--Kung--Morf characterization holds over \emph{every}
commutative ring (\cref{lem:1979paper-lemma1-PIR}):
\[
\alpha_{+}(A)=\rho\bigl(\nabla_{+}(A)\bigr),
\qquad
\alpha_{-}(A)=\rho\bigl(\nabla_{-}(A)\bigr).
\]
Consequently every Toeplitz matrix over any commutative ring has $\alpha_{\pm}\le 2$
(\cref{cor:toeplitz-inverse-ranks}), and the two displacement ranks of any matrix differ by at
most two (\cref{lem:plus-minus-comparison}). We call an elementary divisor ring in which Smith normal
form is unique up to units a \emph{Smith ring}; this includes all principal ideal domains and
principal ideal rings such as $\mathbb{Z}/m\mathbb{Z}$. Over a Smith ring, inner rank equals
the number of nonzero invariant factors (\cref{lem:inner-rank-invariant-factors}), so
displacement ranks can be computed from a Smith normal form.

\smallskip
\noindent\textit{Inversion.}
Over a Smith ring, every invertible $A$ satisfies $\alpha_{+}(A)=\alpha_{-}(A^{-1})$ and
$\alpha_{-}(A)=\alpha_{+}(A^{-1})$ (\cref{lem:displacement-inversion}). The proof reduces the
statement to showing that $I-BC$ and $I-CB$ have the same number of nonzero invariant factors.
We establish this through a block-matrix equivalence together with uniqueness of Smith normal
form, which replaces the rank identity used over fields (\cref{cor:key-equality}). As a
consequence, the inverse of an invertible Toeplitz matrix over a Smith ring, for example over
$\mathbb{Z}$ or $\mathbb{Z}/m\mathbb{Z}$, is again a sum of at most two products of triangular
Toeplitz matrices. Whether inversion exchanges $\alpha_{+}$ and $\alpha_{-}$ over an arbitrary
commutative ring is left open.

\smallskip
\noindent\textit{Hankel and cyclic displacement.}
Right multiplication by the reversal matrix turns the Hankel displacement equation $A-ZAZ=G$ into the
Toeplitz one, which gives an explicit unique solution and identifies the minimal length of the
corresponding representation with $\rho(A-ZAZ)$ over any commutative ring
(\cref{lem:hankel-displacement,cor:hanhle-disp-cor}). The cyclic shift $Z_c=Z+e_1e_n^{\top}$ behaves
differently, because $Z_c$ is invertible rather than nilpotent. The displacement rank under
$Z_c$ differs from $\alpha_{+}$ by at most two
(\cref{lem:comparison-nilpotent-cyclic-displacement}). However, the equation
$A-Z_cAZ_c^{\top}=G$ is solvable if and only if every cyclic diagonal of $G$ sums to zero
(\cref{lem:cyclic-displacement-solvability}), and in that case its solutions are determined
only up to an arbitrary circulant matrix. We write the general solution explicitly as the
Toeplitz reconstruction plus a strictly upper-triangular Toeplitz boundary correction plus a
circulant (\cref{lem:circulant-displacement}). When $n$ is a unit and $R$ contains a principal
$n$-th root of unity, the discrete Fourier transform diagonalizes every circulant matrix, and a
circulant is invertible if and only if all of its Fourier coefficients are units
(\cref{thm:DFT-diagonalization,cor:circulant-inverse}). 

\smallskip
\noindent\textit{Levinson-type recursions.}
For Toeplitz matrices whose leading principal submatrices are all invertible, we prove
Levinson-type recursions over an arbitrary commutative ring. The key point in each case is
that every quantity inverted in the recursion is a unit of $R$, not merely nonzero. For
symmetric block Toeplitz matrices, the residual blocks are Schur complements of consecutive
leading principal submatrices and are therefore invertible
(\cref{lem:levinson-residual-invertibility,thm:generalized-levinson-algorithm}). For
nonsymmetric scalar Toeplitz systems, Jacobi's complementary minor identity gives
\[
d_k=\frac{\det(T_{k+1})\,\det(T_{k-1})}{\det(T_k)^{2}}
\]
for the recursion denominators, so each $d_k$ is a unit (\cref{lem:denominator,thm:nonsym}).
The resulting solver computes $T_n^{-1}b$ using $O(n^{2})$ ring operations, at most $n$
inversions of units, and $O(n)$ working storage (\cref{cor:levinson-computational-cost}). An
example over $\mathbb{Z}/6\mathbb{Z}$ shows that invertibility of $T_n$ alone does not
suffice.

\smallskip
We do not treat nonsymmetric block Toeplitz recursions, superfast algorithms, or bit
complexity, and all operation counts are in ring operations.

\subsection{Related work}
\label{ssec:related-work}

Recursive methods for Toeplitz systems preceded the formulation of displacement rank.
Levinson introduced his recursion in the context of linear prediction and filter design
\cite{Levinson46}, and Trench gave an algorithm for inverting finite Toeplitz matrices
\cite{Trench64}. Bareiss developed elimination methods for Toeplitz and vector Toeplitz
systems \cite{Bareiss69}. These methods exploit Toeplitz structure directly.

Kailath, Kung, and Morf introduced displacement rank as a measure of departure from Toeplitz
structure and related low-rank factorizations of the displacement to sums of products of
triangular Toeplitz matrices \cite{KKM79}. This extended fast methods to matrices that are not
Toeplitz but have small displacement rank. The framework was later extended to Hankel,
Vandermonde, and related matrices arising through products and inverses \cite{CK91}, and to
block-Toeplitz, Toeplitz-block, and Toeplitz-derived matrices through generalized
displacement structure and Schur-type algorithms \cite{KC94}. Kailath and Sayed surveyed these
developments and their connections to generalized Schur algorithms, triangular factorization,
interpolation, and filtering \cite{KS95}. Their treatment of displacement structure under
Schur complementation underlies recursive factorization methods.

Displacement structure also led to algorithms of lower complexity. Bitmead and Anderson gave
an inversion algorithm for generic Toeplitz systems using $O(n\log^{2}n)$ operations
\cite{BA80}, and Morf developed doubling algorithms in which larger structured problems are
solved by combining solutions of smaller ones \cite{Morf80}. Gohberg, Kailath, and Olshevsky used transformations between
displacement structures to obtain fast Gaussian elimination with partial pivoting for
Cauchy-like matrices and, through them, for other structured classes \cite{GKO95}. Bostan,
Jeannerod, and Schost used fast matrix multiplication to improve the dependence of the cost on
the displacement rank \cite{BJS08}. Bostan, Jeannerod, Mouilleron, and Schost then generalized
the displacement operators to block diagonal matrices with companion blocks \cite{BJMS17}. On
the numerical side, Bojanczyk, Brent, de Hoog, and Sweet established stability of the Bareiss
algorithm for symmetric positive definite Toeplitz matrices \cite{BBHS95}. Chandrasekaran, Gu,
Sun, Xia, and Zhu built a superfast solver by passing to Cauchy-like matrices and exploiting
low numerical rank of their off-diagonal blocks \cite{CGSXZ07}.

Most of this literature works over $\mathbb{R}$, $\mathbb{C}$, or an arbitrary field. The
closest precedent for the present setting is Bareiss's study of matrix problems over integral
domains, which includes exact solution of finite Toeplitz systems \cite{Bareiss72}. We work
over general commutative rings, which may have zero divisors. We replace rank by inner rank,
relate it to Smith invariant factors where possible, and make explicit the unit conditions that
the recursions require.

\subsection{Organization}
\label{ssec:organization}

\Cref{sec:prelim} reviews the ring-theoretic background, introduces Smith rings and inner rank,
and defines displacement ranks. \Cref{sec:characterization} proves the inner-rank characterization of
displacement ranks and the inversion and comparison results. \Cref{sec:hankel-cyclic} treats
Hankel and cyclic displacement equations and the Fourier diagonalization of circulant
matrices. \Cref{sec:levinson} develops the Levinson-type recursions and their cost.

\subsection*{Acknowledgments} R.S. was supported by the U.S. Department of Energy, Office of Science, Accelerated Research in Quantum Computing Centers, Quantum Utility through Advanced Computational Quantum Algorithms, grant no. DE-SC0025572.

\subsection*{On AI use}
Generative AI was used to assist in editing portions of the paper, and to check the proofs. The authors reviewed and revised the generated text and take responsibility for the mathematical content, citations, and final manuscript.

\section{Displacement ranks of matrices over rings}
\label{sec:prelim}

This section introduces the background used throughout the paper. In \cref{ssec:comm-rings-ranks}, we review the relevant ring structures and notions of rank, and in \cref{ssec:displacement-rank}, we introduce Toeplitz structure and displacement rank.

\subsection{Commutative rings and matrix ranks}
\label{ssec:comm-rings-ranks}

We begin by recalling some standard facts about commutative rings and matrices over them, restricting our discussion to those aspects necessary for this paper. The interested reader is referred to standard textbooks \cite{Lang02,Brown93} for a more thorough treatment of these topics. Some basic notation used throughout the paper is introduced in this subsection, while additional notation is introduced later as needed. A \emph{commutative ring} $R$ is a ring in which both addition and multiplication are commutative. Throughout this paper, we assume that $R$ has a multiplicative identity, denoted by $1$. For $a,b \in R$, we write $a \mid b$ if there exists $c \in R$ such that $ac = b$. The set of all invertible elements of $R$, called \emph{units}, is denoted by $R^\times$. Thus, $a \in R^\times$ if and only if there exists $b \in R$ such that $ab = 1$. Recall that $R^\times$ forms an abelian group under multiplication, called the \emph{group of units}.

The set of all $m \times n$ matrices over $R$ is denoted $R^{m \times n}$. If $A \in R^{m \times n}$, we denote its transpose as $A^\top$, and it is an element of $R^{n \times m}$. The $(i,j)$ entry of $A$ is denoted $A_{ij}$, or sometimes $A_{i,j}$. The notation $\text{diag} (d_1, d_2,\dots,d_n)$ means the usual shorthand for an $n \times n$ diagonal matrix, whose diagonal entries are $d_1,d_2,\dots,d_n$. Similarly, a matrix $A$ is \emph{block diagonal} if $A=\operatorname{diag}(A_1,\dots,A_k)$, where each $A_i$ is a square matrix and all off-diagonal blocks are zero; here, the notation $\operatorname{diag}(A_1, \dots, A_k)$ is a shorthand for the direct sum notation $A_1 \oplus \cdots \oplus A_k$. We denote the $n \times n$ identity matrix as $I_n$, and when its shape can be inferred from context, we simply write $I$. The set of all invertible $n \times n$ matrices in $R^{n \times n}$ form a group, and is denoted $\operatorname{GL}_n(R)$. If $A \in R^{n \times n}$, we denote its determinant as $\operatorname{det}(A)$, and we recall that $A \in \operatorname{GL}_n(R)$ if and only if $\operatorname{det}(A) \in R^\times$. When $A$ is invertible, its inverse is denoted $A^{-1}$, and it is unique, and satisfies $AA^{-1} = A^{-1}A=I$. A square matrix $N$ is called \emph{nilpotent} if $N^s=0$, for some positive integer $s$.

If $x \in R^m$, its entries will be denoted $x_1,\dots,x_m$, i.e. $x := (x_1,\cdots,x_m)$, and for convenience, we will regard $x$ as an element of $R^{m \times 1} \simeq R^m$; therefore, we treat $x \in R^m$ as a column matrix, and $x^\top \in R^{1\times m} \simeq R^m$ as a row matrix, unless stated otherwise. We denote by $e_i$ the \emph{standard basis column matrix}, whose $i$-th entry is $1$, and all other entries are zero, and its dimension is determined from context. Similarly, for matrices, we define \(E_{ij} \in R^{m \times m}\) (with $m$ determined from context) to be the matrix with \(1\) in position \((i,j)\), and zeros everywhere else.

A \emph{principal ideal ring} (PIR) is a commutative ring in which every ideal is generated by a single element. An \emph{integral domain} is a nonzero commutative ring where there are no zero divisors. A \emph{principal ideal domain} (PID) is a PIR that is also an integral domain. An \emph{elementary divisor ring} (EDR) is a commutative ring in which every matrix has a Smith normal form, as defined below: 
\begin{defn}[Smith normal form]
\label{def:snf}
Let $R$ be an elementary divisor ring, and let $A \in R^{m \times n}$. Then, a \emph{Smith normal form} of $A$ is a diagonal matrix $D=\operatorname{diag}(d_1, \ldots , d_r, 0, \ldots , 0) \in R^{m \times n}$, with each $d_i \neq 0$, such that
$A= PDQ$, for some invertible matrices $P \in \operatorname{GL}_m(R)$ and $Q \in \operatorname{GL}_n(R)$, and where $d_1 \mid d_2 \mid \cdots \mid d_r$. The nonzero diagonal entries $d_1, \ldots, d_r$ are called the \emph{invariant factors} of $A$, and we define $\mu(A) := r \le \min\{m,n\}$, the number of nonzero invariant factors (this is well defined over a Smith ring by the uniqueness property stated below). The factorization $A=PDQ$ is called a \emph{Smith normal form decomposition} of $A$.
\end{defn}

Several arguments below use the following uniqueness property. If $A=PDQ=P'D'Q'$ are two Smith normal form decompositions of $A\in R^{m\times n}$ with $D=\operatorname{diag}(d_1,\dots,d_r,0,\dots,0)$ and
$D'=\operatorname{diag}(d_1',\dots,d_s',0,\dots,0)$, then $r=s$ and $d_i=u_id_i'$ for some units $u_i\in R^\times$, $1\le i\le r$. This is classical for principal ideal domains and holds for principal ideal rings such as $\mathbb Z_d$ \cite{Brown93}. For a general elementary divisor ring we do not claim it, because when $R$ has zero divisors the number of nonzero invariant factors cannot be read off from determinantal ideals. We call an elementary divisor ring in which this uniqueness holds a \emph{Smith ring}. This
includes principal ideal domains and the principal ideal rings mentioned above. Whenever Smith uniqueness is used below, we assume that $R$ is a Smith ring.

It is known that every PIR is also an EDR, but the converse is not true --- for example, the commutative ring of entire holomorphic functions is an EDR, but not a PIR. Examples of PIR include the quotient ring $\mathbb Z_d := \mathbb Z /d \mathbb Z$, while the ring of integers $\mathbb Z$ is an example of a PID. Over a general commutative ring, there are several different notions of rank of a matrix. One such notion of rank, which will be useful to us in proofs, is introduced below.

\begin{defn}[Inner rank]
\label{def:inner-rank}
Let $R$ be a commutative ring, and $A\in R^{m\times n}$. The \emph{inner rank} of $A$, denoted $\rho(A)$, is the least integer $k\ge 0$, such that there exist matrices 
$B\in R^{m\times k}$, and $C\in R^{k\times n}$, such that $A = B\,C$. We adopt the convention $\rho(0)=0$.
\end{defn}

It turns out that for a Smith ring, the notion of inner rank of a matrix coincides with the number of invariant factors, and this is established in the following lemma:

\begin{lemma}
\label{lem:inner-rank-invariant-factors}
Let $R$ be a Smith ring, and $A\in R^{m\times n}$. Let $P_A\in\mathrm{GL}_m(R)$ and $Q_A\in\mathrm{GL}_n(R)$ be such that
$A=P_A D_A Q_A$ is a Smith normal form of $A$, where
$D_A :=\operatorname{diag}(d_1,\dots,d_r,0,\dots,0)$,
with $d_i\neq 0$, and $d_i\mid d_{i+1}$, for
$1\leq i<r$. Then, $\rho(A)=\mu(A) = r$.
\end{lemma}

\begin{proof}
If $r=0$, then $D_A=0$, so $A=0$ and $\rho(A)=0$, so we may assume that $r>0$, in the rest of the proof. Let $E\in R^{m\times r}$ be the matrix consisting of the first $r$ columns of $D_A$, and let $F\in R^{r\times n}$ be the matrix consisting of the first $r$ rows of $I_n$. Then, $D_A=EF$, so $A=(P_AE)(FQ_A)$. This shows that $\rho(A)\leq r$, and thus, we only need to show that $\rho(A)\ge r$.

Let us assume for contradiction, that $A=BC$, with $B\in R^{m\times k}$ and $C\in R^{k\times n}$, for some $k<r\leq\min\{m,n\}$. Since $A\neq 0$ (as $r > 0$), we must have $k\geq 1$. By the Smith normal form decomposition, we can assume that $B=P_BD_BQ_B$ and $C=P_CD_CQ_C$, for some $P_B\in\mathrm{GL}_m(R)$, $Q_B,P_C\in\mathrm{GL}_k(R)$, $Q_C\in\mathrm{GL}_n(R)$, and diagonal matrices $D_B\in R^{m\times k}$, and $D_C\in R^{k\times n}$. Then, $A=P_BXQ_C$, where $X:=D_B (Q_BP_C) D_C \in R^{m \times n}$. Since $D_B$ and $D_C$ are rectangular diagonal matrices, the last $m-k$ rows of $D_B$, and the last $n-k$ columns of $D_C$ are zero. Hence, we have the block diagonal form $X=\left( \begin{smallmatrix}
V & 0\\
0 & 0
\end{smallmatrix} \right)$, where $V\in R^{k\times k}$. Now, take a Smith normal form $V=P_VD_VQ_V$, where
$P_V,Q_V\in\mathrm{GL}_k(R)$, and a diagonal matrix $D_V:=\operatorname{diag}(c_1,\dots,c_s,0,\dots,0)$, with $c_i\neq 0$, $c_i\mid c_{i+1}$ for
$1\leq i<s$, and $s\leq k$. If we extend $P_V^{-1}$ and $Q_V^{-1}$ by identity blocks, we may define $\widehat P:=\operatorname{diag}(P_V^{-1},I_{m-k})$, and $\widehat Q:=\operatorname{diag}(Q_V^{-1},I_{n-k})$. Note that both $\widehat P$ and $\widehat Q$ are also invertible, since they are block diagonal, and their individual blocks are invertible. This leads us to the identity $\widehat P X\widehat Q =\operatorname{diag}(c_1,\dots,c_s,0,\dots,0) =:D' \in R^{m \times n}$. Since we know that $X=P_B^{-1}AQ_C^{-1}$, it follows that
$(\widehat P P_B^{-1})A(Q_C^{-1}\widehat Q)=D'$.
With both left and right multiplying matrices being invertible, this implies that $D'$ is a Smith normal form of $A$. Since the Smith normal form is unique up to multiplication of the diagonal entries by units, we know $D'$ and $D_A$ have the same number of nonzero entries. Thus, $r=s\leq k$, which is a contradiction.
\end{proof}

\subsection{Displacement rank}
\label{ssec:displacement-rank}

Let $R$ be a commutative ring. Then $A \in R^{m \times n}$ is called a \emph{Toeplitz matrix}  if its entries are constant along the diagonal, and every sub- and super-diagonals, that is, $A_{ij} = A_{i+1, \; j+1}$, for all $1 \le i < m$, and $1 \le j < n$. If $x, y \in R^m$, we will use the notation \(L(x)\) to denote a $m \times m$ lower-triangular  Toeplitz matrix, determined entirely by its first column \(x\), while the notation \(U(y^\top)\) will be used to denote a $m \times m$ upper-triangular Toeplitz matrix, determined entirely by its first row \(y^\top\).

\begin{defn}[Displacement rank]
\label{def:displacement-rank}
Let $R$ be a commutative ring, and let $A\in R^{n \times n}$. The \emph{$(+)$-displacement rank}, written $\alpha_+(A)$, and the \emph{$(-)$-displacement rank}, written $\alpha_-(A)$, are defined to be the smallest nonnegative integers, for which $A$ has the representations
\[
A=\sum_{i=1}^{\alpha_+(A)}L_iU_i,
\qquad\text{and}\qquad
A=\sum_{i=1}^{\alpha_-(A)}\widehat U_i\widehat L_i,
\]
respectively, where $L_i,\widehat L_i\in R^{n \times n}$ are lower-triangular Toeplitz matrices, and $U_i,\widehat U_i\in R^{n \times n}$ are upper-triangular
Toeplitz matrices. The case, where the sums are empty, corresponds to the zero matrix, in which case we set $\alpha_+(0)=\alpha_-(0)=0$.
\end{defn}

In \cite{KKM79}, Kailath, Kung, and Morf proved the following theorem about displacement rank of a matrix, in the setting $\mathbb{K} = \mathbb R$ or $\mathbb C$, although most of their results generalize to the setting where $\mathbb{K}$ is a field. For the statement of the theorem, we need to define the \emph{lower-shift matrix} $Z\in R^{n\times n}$ as $Z := \sum_{i=1}^{n-1} E_{i+1,i}$. 

\begin{theorem}[Kailath-Kung-Morf]
\label{thm:disp-rank-fields}
Let $\mathbb{K}=\mathbb{R}$ or $\mathbb{C}$, and let
$A\in\mathbb{K}^{n\times n}$. Let $Z$ be the
$n\times n$ lower-shift matrix. Then we have the characterizations $\alpha_+(A)=\operatorname{rank}(A-ZAZ^\top)$ and $\alpha_-(A)=\operatorname{rank}(A-Z^\top AZ)$.
Furthermore, if $A$ is invertible, then $\alpha_+(A)=\alpha_-(A^{-1})$, and $\alpha_-(A)=\alpha_+(A^{-1})$. 
\end{theorem}

Thus, a matrix with small displacement ranks retains this property under inversion, even though its inverse need not be Toeplitz. In particular, if $A$ is an invertible
Toeplitz matrix, then both $A$ and $A^{-1}$ have displacement ranks at most two $\alpha_+(A),\ \alpha_-(A),\ \alpha_+(A^{-1}),\ \alpha_-(A^{-1})\leq 2$.








\section{Displacement Ranks as Inner Ranks}
\label{sec:characterization}

\subsection{The Inner-Rank Characterization}

This section characterizes displacement ranks through
structured matrix representations and inner rank.
We then use Smith normal form identities to study
inversion and compare the positive and negative
displacement ranks.

\begin{lemma}
\label{lem:inner-rank-properties}
Let $R$ be a commutative ring with identity. For matrices over $R$, inner rank satisfies the following properties:
\begin{enumerate}
    \item $\rho(PAQ)=\rho(A)$ whenever $P$ and $Q$
    are invertible matrices of the appropriate sizes.
    \item $\rho(BAC)\leq\rho(A)$ whenever the product
    $BAC$ is defined.
    \item $\rho(A+B)\leq\rho(A)+\rho(B)$ whenever
    $A$ and $B$ have the same dimensions.
\end{enumerate}
\end{lemma}

\begin{proof}
Let $A=XY$ be a factorization with $X\in R^{m\times r}$ and $Y\in R^{r\times n}$, where $r=\rho(A)$. Then $BAC=(BX)(YC)$, so $\rho(BAC)\leq r=\rho(A)$, proving the second property.

Now, let us start with $\rho(PAQ)\leq\rho(A)$. If $P$ and $Q$ are invertible, applying the same inequality to $A=P^{-1}(PAQ)Q^{-1}$ results in $\rho(A)\leq\rho(PAQ)$. Thus, $\rho(PAQ)=\rho(A)$, proving the first property.

Finally, let $A=XY$ and $B=UV$ be factorizations through $R^r$ and $R^s$, respectively, such that $r=\rho(A)$ and $s=\rho(B)$. Then $A+B=\begin{bmatrix}X&U\end{bmatrix}
\left[\begin{smallmatrix}Y\\V\end{smallmatrix}\right]$.
Since this is a factorization through $R^{r+s}$, we know $\rho(A+B)\leq r+s=\rho(A)+\rho(B)$, proving the last property.
\end{proof}

\begin{lemma}
\label{lem:nilpotent-displacement}
Let ${R}$ be a commutative ring with identity.
Let $P\in{R}^{m\times m}$,
$Q\in{R}^{n\times n}$, and
$G\in{R}^{m\times n}$.
Suppose that $P^s=0$ or $Q^s=0$ for some integer
$s\geq 1$.
Then the equation $X-PXQ=G$ has the unique solution
\[
X=\sum_{k=1}^{s}P^{k-1}GQ^{k-1}
\]
\end{lemma}

\begin{proof}
To prove uniqueness, suppose that $X_1$ and $X_2$
are solutions, and let $\Delta=X_1-X_2$.
Subtracting their equations gives $\Delta=P\Delta Q$.
By induction, $\Delta=P^k\Delta Q^k$ for every
integer $k\geq 1$.
Taking $k=s$, we obtain $\Delta=0$, since
$P^s=0$ or $Q^s=0$.
Hence $X_1=X_2$.

To verify the proposed solution, let $X=\sum_{k=1}^{s}P^{k-1}GQ^{k-1}$. Then we will use the following sequence of equalities $X-PXQ =\sum_{k=1}^{s}P^{k-1}GQ^{k-1}-\sum_{k=2}^{s+1}P^{k-1}GQ^{k-1} = G-P^sGQ^s = G$, where the last equality holds because $P^s=0$ or $Q^s=0$. Thus $X$ is a solution, and by the first part it is the unique one.
\end{proof}

For the remainder of this section, we will write $\nabla_+(A)=A-ZAZ^\top$ and $\nabla_-(A)=A-Z^\top AZ$, where $Z$ is the standard lower-shift matrix.

\begin{lemma}
\label{lem:unique-sol-pir} 
Let $R$ be a commutative ring with identity, $\{x_i,y_i\in R^n:i=1,\dots,\alpha\}$ be column vectors, and let $Z\in R^{n\times n}$ be the standard lower-shift matrix. Then, $A - Z A Z^{\top} \;=\; \sum_{i=1}^{\alpha} x_i y_i^{\top}$ has the unique solution
\begin{equation}
\label{eq:unique-sol-pir-1}
  A \;=\; \sum_{i=1}^{\alpha} L(x_i)\,U(y_i^{\top}),
\end{equation}
where $L(x)$ is the lower–triangular Toeplitz matrix over $R$, whose first column is $x$, and $U(y^{\top})$ is the upper–triangular Toeplitz matrix over $R$, whose first row is $y^{\top}$.
Similarly, $A - Z^{\top} A Z \;=\; \sum_{i=1}^{\alpha} x_i y_i^{\top}$ has the unique solution
\begin{equation}
\label{eq:unique-sol-pir-2}
  A \;=\; \sum_{i=1}^{\alpha} U  \big({\tilde{x_i}}^{\top}\big)\,L \big({\tilde{y_i}}\big),
\end{equation}
where ${\tilde{x}} = [x_n\,\, x_{n-1}\,\,\cdots\,\, x_1]^{\top}$ when $x = [x_1\,\, x_2\,\,\cdots\,\, x_n]^{\top}$.
\end{lemma}

\begin{proof}
Uniqueness follows from Lemma \ref{lem:nilpotent-displacement}, taking $P=Z$ and $Q=Z^\top$ for the first equation, and $P=Z^\top$ and $Q=Z$ for the second equation, since $Z^n=(Z^\top)^n=0$.

We find the solution of $A-ZAZ^\top=\sum_{i=1}^{\alpha}x_i y_i^\top$ by showing $L(x)U(y^\top) -Z[L(x)U(y^\top)]Z^\top=xy^\top$, to verify \eqref{eq:unique-sol-pir-1}. Since $ZL(x)=L(x)Z$, $U(y^\top)Z^\top=Z^\top U(y^\top)$, and $I-ZZ^\top=e_1e_1^\top$, we know $L(x)U(y^\top) -Z[L(x)U(y^\top)]Z^\top = xy^\top$. When we sum this identity over $i=1,\dots,\alpha$, it shows that \eqref{eq:unique-sol-pir-1} satisfies the equation.

Now, in order to verify \eqref{eq:unique-sol-pir-2}, we will use $U=U(\tilde{x}^\top)$ and $L=L(\tilde{y})$. By the definitions of the Toeplitz matrices, their last column and last row should satisfy $Ue_n=x$ and $e_n^\top L=y^\top$. Also, since $Z^\top U=UZ^\top$, $LZ=ZL$, and $I-Z^\top Z=e_ne_n^\top$, we derive $UL-Z^\top ULZ =U(I-Z^\top Z)L = Ue_ne_n^\top L = xy^\top$. When we sum over $i=1,\dots,\alpha$, it verifies \eqref{eq:unique-sol-pir-2}, as desired.
\end{proof}

\begin{lemma}
\label{lem:1979paper-lemma1-PIR}
Let $R$ be a commutative ring and let $A\in {R}^{n\times n}$.
Let $Z\in{R}^{n\times n}$ be the lower-shift.
Then the (+)-displacement rank can be computed by
\begin{equation}
\label{plus-displacement-rank}
  \alpha_+(A)=\rho(A-ZAZ^\top)
\end{equation}
and the (-)-displacement rank can be computed by
\begin{equation}
\label{minus-displacement-rank}
 \alpha_-(A)=\rho(A-Z^\top AZ)
\end{equation}
where $\rho$ denotes inner rank, as in Definition~\ref{def:inner-rank}.
\end{lemma}

\begin{proof}
Assume $A$ has a minimal $(+)$-representation,
$A=\sum_{i=1}^{\alpha_{+}(A)} L(x_i)\,U(y_i^\top)$.
By the standard Toeplitz shift identities, $ZL(x_i)=L(x_i)Z$ and $U(y_i^\top)Z^\top=Z^\top U(y_i^\top)$. Since $I-ZZ^\top=e_1e_1^\top$, we know $L(x_i)U(y_i^\top)-ZL(x_i)U(y_i^\top)Z^\top = L(x_i)(I-ZZ^\top)U(y_i^\top) = x_i y_i^\top.$
Hence, we know that $A-ZAZ^\top=\sum_{i=1}^{\alpha_{+}(A)}x_i y_i^\top$.
Writing these outer products as a matrix factorization gives us
\[
A-ZAZ^\top
=
\begin{bmatrix}x_1&\cdots&x_{\alpha_+(A)}\end{bmatrix}
\begin{bmatrix}
y_1^\top\\
\vdots\\
y_{\alpha_+(A)}^\top
\end{bmatrix}.
\]
Then, by the definition of inner rank, we know that $\rho(A-ZAZ^\top)\leq\alpha_+(A)$.

If $\rho(A-ZAZ^\top)=r$, then there exist
$x_i,y_i\in R^n$ such that $A-ZAZ^\top=\sum_{i=1}^{r}x_i y_i^\top$.
Using Lemma~\ref{lem:unique-sol-pir}, we know $A=\sum_{i=1}^{r}L(x_i)\,U(y_i^\top)$,
so $\alpha_+(A)\leq r$. Combined with the first inequality, this proves the $(+)$ equality in \eqref{plus-displacement-rank}.

The proof for the $(-)$-displacement rank is analogous,
with $Z$ and $Z^\top$ interchanged.
For the upper-triangular Toeplitz matrix $U$ and the lower-triangular Toeplitz matrix $L$, we derive that $UL-Z^\top ULZ = U(I-Z^\top Z)L = (Ue_n)(e_n^\top L)$, which is an outer product. Thus, a minimal $(-)$-representation gives
$\rho(A-Z^\top AZ)\leq\alpha_-(A)$.

The negative case of Lemma~\ref{lem:unique-sol-pir} converts a factorization of $A-Z^\top AZ$ of inner rank $r$ into a sum of $r$ upper-triangular Toeplitz times
lower-triangular Toeplitz products. Hence
$\alpha_-(A)\leq r$, which together with the first inequality proves the $(-)$ equality in \ref{minus-displacement-rank}.
\end{proof}

\subsection{Inversion and Comparison of the Two Displacement Ranks}

\begin{lemma}
\label{lem:key-equality-invariant-factors}
Let $R$ be a Smith ring. Let $A \in R^{m \times n}$ and $B \in R^{n \times m}$ be two matrices. Then the two block matrices
\[
    \begin{pmatrix}
        I_{n} & 0 \\
        0 & I_m - AB
    \end{pmatrix}
    \quad \text{and} \quad
    \begin{pmatrix}
        I_{m} & 0 \\ 
        0 & I_n - BA
    \end{pmatrix}
\]
have the same Smith normal form, up to multiplication by units.
\end{lemma}

\begin{proof}
Consider the $(m+n)\times(m+n)$ block matrix $M=\begin{pmatrix}
I_m & A\\
B & I_n
\end{pmatrix}$.
We want to show that $M$ is equivalent to both block
matrices given above by performing row and column
operations. Here, equivalence means multiplication
on the left and right by invertible matrices.

We start by eliminating the block $B$ in the lower-left corner, by left-multiplying $M$ with an invertible
lower-triangular block matrix as follows
\[
\begin{pmatrix}
I_m & 0\\
-B & I_n
\end{pmatrix}
\begin{pmatrix}
I_m & A\\
B & I_n
\end{pmatrix}
=
\begin{pmatrix}
I_m & A\\
0 & I_n-BA
\end{pmatrix}.
\]

Now, we can eliminate the block $A$ in the upper-right
corner by right-multiplying by an invertible
upper-triangular block matrix. This gives
\[
\begin{pmatrix}
I_m & A\\
0 & I_n-BA
\end{pmatrix}
\begin{pmatrix}
I_m & -A\\
0 & I_n
\end{pmatrix}
=
\begin{pmatrix}
I_m & 0\\
0 & I_n-BA
\end{pmatrix}.
\]
Hence, $M\sim\operatorname{diag}(I_m,I_n-BA)$.

Separately, from $M$, we can eliminate the block $A$
in the upper-right corner by left-multiplying
\[
\begin{pmatrix}
I_m & -A\\
0 & I_n
\end{pmatrix}
\begin{pmatrix}
I_m & A\\
B & I_n
\end{pmatrix}
=
\begin{pmatrix}
I_m-AB & 0\\
B & I_n
\end{pmatrix}.
\]
Then, we can eliminate the block $B$ in the lower-left
corner by right-multiplication. This gives
\[
\begin{pmatrix}
I_m-AB & 0\\
B & I_n
\end{pmatrix}
\begin{pmatrix}
I_m & 0\\
-B & I_n
\end{pmatrix}
=
\begin{pmatrix}
I_m-AB & 0\\
0 & I_n
\end{pmatrix}.
\]
Since the multiplying matrices used above are invertible, their inverses are found by replacing $-A$ with $A$ and $-B$ with $B$. Thus, $M\sim\operatorname{diag}(I_m-AB,I_n)$.

By permuting the rows and columns to exchange the
diagonal blocks, we also know that $\operatorname{diag}(I_m-AB,I_n) \sim\operatorname{diag}(I_n,I_m-AB)$. Since $\operatorname{diag}(I_m,I_n-BA)$ and
$\operatorname{diag}(I_n,I_m-AB)$ are equivalent to $M$,
they are equivalent to each other. Over a Smith ring, two equivalent matrices have the same Smith normal form up to multiplication of the diagonal entries by units, as desired.
\end{proof}

\begin{cor}
\label{cor:key-equality}
Let $R$ be a Smith ring. Let $A,B \in R^{n \times n}$ be matrices. Then, $I-AB$ and $I-BA$ will have the same Smith normal form, up to multiplication of the diagonal entries by units. In particular, they have the same number of nonzero invariant factors.
\end{cor}

\begin{proof}
By Lemma~\ref{lem:key-equality-invariant-factors},
the matrices $\operatorname{diag}(I_n,I_n-AB) \quad\text{and}\quad
\operatorname{diag}(I_n,I_n-BA)$ have the same Smith normal form, up to multiplication of the diagonal entries by units.

Let $D=\operatorname{diag}(d_1,\dots,d_n)$ and $E=\operatorname{diag}(e_1,\dots,e_n)$ be Smith normal forms of $I_n-AB$ and $I_n-BA$, respectively, including any trailing zero entries. Then, $\operatorname{diag}(I_n,D)$ and $\operatorname{diag}(I_n,E)$ are the Smith normal forms
of the two block matrices above, since $1$ divides
every element of $R$.

By uniqueness of Smith normal form up to multiplication
of the diagonal entries by units, the corresponding
entries agree up to units. Comparing the last $n$
entries, we get $d_i=u_i e_i$ where $i=1,\dots,n$ for some units $u_i\in R^\times$. Then, $I_n-AB$ and $I_n-BA$ have the same Smith normal form up to units. Also, since multiplication by a unit
preserves whether an entry is zero, they also have the same number of nonzero invariant factors.
\end{proof}

\begin{lemma}
\label{lem:displacement-inversion}
Let $R$ be a Smith ring, and let $A\in\mathrm{GL}_n(R)$. Then, $\alpha_+(A)=\alpha_-(A^{-1})$ and also $\alpha_-(A)=\alpha_+(A^{-1})$.
\end{lemma}

\begin{proof}
Let $Z$ be the lower-shift matrix. By the characterization of displacement ranks in Lemma \ref{lem:1979paper-lemma1-PIR}, we know $\alpha_+(A)=\rho(A-ZAZ^\top)$ and $\alpha_-(A^{-1})=\rho(A^{-1}-Z^\top A^{-1}Z)$.

Multiplication on either side by an invertible matrix preserves inner rank. Any factorization of a matrix gives a factorization of its product with an invertible matrix, and the reverse inequality follows by multiplying by the inverse. Hence, $\alpha_+(A) = \rho\bigl(A^{-1}(A-ZAZ^\top)\bigr) = \rho(I_n-A^{-1}ZAZ^\top)$, and similarly $\alpha_-(A^{-1}) = \rho\bigl((A^{-1}-Z^\top A^{-1}Z)A\bigr) = \rho(I_n-Z^\top A^{-1}ZA)$.

Now, set $B=A^{-1}ZA$ and $C=Z^\top$. The two matrices on the right are $I_n-BC$ and $I_n-CB$. Corollary~\ref{cor:key-equality} shows that these matrices will have the same Smith normal form up to units. In particular, they have the same number of nonzero invariant factors. Since we know that  inner rank equals the number of nonzero invariant factors over a Smith ring, their inner ranks are equal. Hence, $\alpha_+(A)=\alpha_-(A^{-1})$.

Applying the equality to $A^{-1}$ results in $\alpha_+(A^{-1})=\alpha_-(A)$, proving the second equality.
\end{proof}

\begin{cor}
\label{cor:toeplitz-inverse-ranks}
Let $R$ be a commutative ring, and let $T\in R^{n\times n}$ be a Toeplitz matrix. Then $\alpha_+(T)\le 2$ and $\alpha_-(T)\le 2$. Moreover, if $R$ is a Smith ring and $T$ is invertible, then $\alpha_+(T^{-1})\le 2$ and $\alpha_-(T^{-1})\le 2$.
\end{cor}

\begin{proof}
Since $T$ is Toeplitz, the entries of $D_+=\nabla_+(T)$ satisfy $(D_+)_{ij}=T_{ij}-T_{i-1,j-1}=0$ when $2\leq i,j\leq n$. Thus, $D_+$ is supported on its first row and column, and we can write $D_+ = e_1(e_1^\top D_+) +\bigl((I-e_1e_1^\top)D_+e_1\bigr)e_1^\top$.
Each term is an outer product and has inner rank
at most one. So, by Lemma~\ref{lem:inner-rank-properties}, we get $\rho(D_+)\leq 2$.

Similarly, the entries of $D_-=\nabla_-(T)$ satisfy $(D_-)_{ij}=T_{ij}-T_{i+1,j+1}=0$ when $1\leq i,j\leq n-1$. Hence, $D_-$ is supported on its last row and column, and $D_- = e_n(e_n^\top D_-) +\bigl((I-e_ne_n^\top)D_-e_n\bigr)e_n^\top$. Therefore, $\rho(D_-)\leq 2$.
By Lemma \ref{lem:1979paper-lemma1-PIR}, we know that $\alpha_+(T)\leq 2$ and $\alpha_-(T)\leq 2$.

If $R$ is a Smith ring and $T$ is invertible, Lemma~\ref{lem:displacement-inversion} tells us
$\alpha_+(T^{-1})=\alpha_-(T)\leq 2$ and $\alpha_-(T^{-1})=\alpha_+(T)\leq 2$, as desired.
\end{proof}

Lemma \ref{lem:displacement-inversion} and the second part of Corollary \ref{cor:toeplitz-inverse-ranks} are proved only over Smith rings. We do not address whether inversion exchanges $\alpha_+$ and $\alpha_-$ over an arbitrary commutative ring.

\begin{lemma}
\label{lem:plus-minus-comparison}
Let $R$ be a commutative ring, and let $A\in R^{n\times n}$. Then
\begin{equation*}
    \abs{\alpha_+(A)-\alpha_-(A)} \le 2.
\end{equation*}
If, in addition, $R$ is a Smith ring, then $\abs{\mu(\nabla_+(A))-\mu(\nabla_-(A))}\le 2$.
\end{lemma}

\begin{proof}
By \cref{lem:1979paper-lemma1-PIR}, $\alpha_\pm(A)=\rho(\nabla_\pm(A))$, so it suffices to compare the inner ranks of the two displacement matrices.

For the standard lower-shift matrix $Z$, we know $ZZ^\top=I-e_1e_1^\top$ and $Z^\top Z=I-e_ne_n^\top$, where $e_1$ and $e_n$
are the first and last standard basis vectors. So, we can express $\nabla_+(A)$ in terms of $\nabla_-(A)$ through $\nabla_+(A)+Z\nabla_-(A)Z^\top = e_1e_1^\top A+Ae_1e_1^\top -e_1e_1^\top Ae_1e_1^\top$.
If we put the boundary terms into two outer products, we get $\nabla_+(A) = -Z\nabla_-(A)Z^\top +e_1(e_1^\top A) +\bigl((I-e_1e_1^\top)Ae_1\bigr)e_1^\top$.

We know that inner rank is subadditive. This means that if $M=BC$ and $N=DE$, then we have $M+N=\begin{bmatrix}B&D\end{bmatrix}
\left[\begin{smallmatrix}C\\E\end{smallmatrix}\right]$. Also, multiplication on either side by a matrix cannot increase inner rank, and each outer product has inner rank at most one. This gives us the inequality $\rho(\nabla_+(A)) \leq \rho(\nabla_-(A))+2$. 

Similarly, we can derive the opposing statement by interchanging $Z$ and $Z^\top$. This will give us $\nabla_-(A)+Z^\top\nabla_+(A)Z = A-Z^\top ZAZ^\top Z = e_n(e_n^\top A) +\bigl((I-e_ne_n^\top)Ae_n\bigr)e_n^\top$. Then, applying the same inner-rank inequalities gives $\rho(\nabla_-(A))\leq\rho(\nabla_+(A))+2$.

Combining the two inequalities gives $\abs{\rho(\nabla_+(A))-\rho(\nabla_-(A))}\le 2$, which is the first claim. If $R$ is a Smith ring, then $\mu=\rho$ by \cref{lem:inner-rank-invariant-factors}, which gives the second claim.
\end{proof}

\section{Hankel and Cyclic Displacement Structure}
\label{sec:hankel-cyclic}

The reconstruction formulas in \cref{sec:characterization} provide a starting point for studying displacement equations associated with other shifts. We first show that multiplication by the reversal matrix converts the Hankel displacement equation into the Toeplitz displacement equation. We then replace the nilpotent lower shift $Z$ by the cyclic shift $Z_c=Z+e_1e_n^\top$. This change introduces a boundary error that can be corrected
using an upper-triangular Toeplitz matrix. Unlike the nilpotent case, the cyclic displacement equation requires a compatibility condition and determines its solution only up to an arbitrary circulant matrix.

\subsection{Hankel Displacement via Reversal}

A matrix $H\in{R}^{n\times n}$ is called a \emph{Hankel matrix} if its entries are constant along each anti-diagonal, so that $H_{ij}=h_{i+j}$ for some elements $h_2,\dots,h_{2n}\in{R}$. If $J$ is the reversal matrix, then $HJ$ is Toeplitz, since $(HJ)_{ij}=h_{n+1+i-j}$ depends only on $i-j$.

\begin{lemma}
\label{lem:hankel-displacement}
Let ${R}$ be a commutative ring with identity and let $Z\in{R}^{n\times n}$ be the standard lower-shift matrix. Let $\{x_i,y_i\in{R}^n:i=1,\dots,\alpha\}$ be column vectors. Then the displacement equation $A-ZAZ=\sum_{i=1}^{\alpha}x_i y_i^\top$ has the unique solution
\[
A=\sum_{i=1}^{\alpha}L(x_i)U((Jy_i)^\top)J,
\]
where $L(x_i)$ is the lower-triangular Toeplitz matrix with first column $x_i$, $U((Jy_i)^\top)$ is the upper-triangular Toeplitz matrix with first row $(Jy_i)^\top$, and $J$ is the reversal matrix, defined by $J_{ij}=1$ if $i+j=n+1$ and $J_{ij}=0$ otherwise.
\end{lemma}

\begin{proof}
Uniqueness follows from Lemma \ref{lem:nilpotent-displacement} with $P=Q=Z$,
since $Z^n=0$.

To obtain the structural representation, we multiply the displacement equation on the right by $J$, which results in $AJ-ZAZJ=\sum_{i=1}^{\alpha}x_i y_i^\top J$. The reversal matrix satisfies $J^\top=J$, $J^2=I$, and $ZJ=JZ^\top$. So, $ZAZJ=ZA(JZ^\top)=Z(AJ)Z^\top$ and $y_i^\top J=(Jy_i)^\top$. Setting $S=AJ$, we obtain
\[
S-ZSZ^\top
=\sum_{i=1}^{\alpha}x_i(Jy_i)^\top.
\]

This is the standard Toeplitz displacement equation. By Lemma \ref{lem:unique-sol-pir}, its unique solution is
\[
S=\sum_{i=1}^{\alpha}L(x_i)U((Jy_i)^\top).
\]
Since $A=SJ$, this gives us the proposed formula.

To verify this formula satisfies the original equation, note that $JZ=Z^\top J$. Then, this gives $A-ZAZ = SJ-ZSJZ = SJ-ZSZ^\top J =(S-ZSZ^\top)J = \sum_{i=1}^{\alpha}x_i(Jy_i)^\top J = \sum_{i=1}^{\alpha}x_i y_i^\top$, where the last equality uses $J^\top J=I$. Thus, the proposed matrix is the unique solution.
\end{proof}

\begin{cor}
Let ${R}$ be a commutative ring with identity, and let $H\in{R}^{n\times n}$ be a Hankel matrix. Then $\rho(H-ZHZ)\leq 2$ and there exist vectors $x_i,y_i\in{R}^n$, for $i=1,2$, such that
\[
H=\sum_{i=1}^{2}L(x_i)U((Jy_i)^\top)J.
\]
\end{cor}

\begin{proof}
Let $G=H-ZHZ$. For $2\leq i\leq n$ and $1\leq j\leq n-1$, we have
$G_{ij}=H_{ij}-H_{i-1,j+1}=0$, since $H$ is Hankel. Thus, $G$ is supported on its first row and last column, and $G=e_1(e_1^\top G) +\bigl((I-e_1e_1^\top)Ge_n\bigr)e_n^\top$. This expresses $G$ as a sum of two outer products, so $\rho(G)\leq 2$. Applying the preceding lemma
to this decomposition gives the stated representation.
\end{proof}

\begin{cor}
\label{cor:hanhle-disp-cor}
Let $R$ be a commutative ring and $A\in R^{n\times n}$. The smallest $r$ for which it holds that 
$A=\sum_{i=1}^{r}L(x_i)U((Jy_i)^\top)J$ for some $x_i,y_i\in R^n$ equals $\rho(A-ZAZ)$.
\end{cor}

\begin{proof}
If $A$ has such a representation with $r$ terms, then we know by Lemma \ref{lem:hankel-displacement} that it satisfies $A-ZAZ=\sum_{i=1}^r x_iy_i^\top$. Then, $\rho(A-ZAZ)\le r$. Conversely, if $\rho(A-ZAZ)=r$, we can write $A-ZAZ=\sum_{i=1}^r x_iy_i^\top$, so Lemma \ref{lem:hankel-displacement} gives a representation with $r$ terms.
\end{proof}

\subsection{Cyclic Displacement: Solvability and the General Solution}
\label{ss:4.2}

The \emph{lower cyclic shift matrix} $Z_c\in{R}^{n\times n}$ is defined by
$Z_ce_j=e_{j+1}$ for $1\leq j<n$ and $Z_ce_n=e_1$.
Equivalently, $Z_c=Z+e_1e_n^\top$, where $Z$ is
the standard nilpotent lower-shift matrix. A matrix $C\in{R}^{n\times n}$ is \emph{circulant}
if $C=\sum_{k=1}^{n}c_kZ_c^{k-1}$ for some $c_1,\dots,c_n\in{R}$. Its first column is $(c_1,\dots,c_n)^\top$. For an $n\times n$ matrix, the \emph{cyclic diagonal} with offset $d\in\{0,\dots,n-1\}$ consists of the positions $(i,j)$ satisfying $j-i\equiv d\pmod n$, where $1\leq i,j\leq n$. A circulant matrix is constant along each cyclic diagonal.

\begin{lemma}
\label{lem:comparison-nilpotent-cyclic-displacement}
Let ${R}$ be a commutative ring, and let $A\in{R}^{n\times n}$. Let $Z$ be the standard lower-shift matrix and $Z_c = Z + e_1 e_n^\top$ be the lower cyclic shift matrix. Define $\alpha_c(A) = \rho(A - Z_c A Z_c^\top)$ as the displacement rank under $Z_c$. 
Then
\[
\alpha_+(A) - 2 \le \alpha_c(A) \le \alpha_+(A) + 2
\]
\end{lemma}

\begin{proof}
Write $Z_c=Z+Q$, where $Q=e_1e_n^\top$. Substituting into the cyclic displacement equation gives $A-Z_cAZ_c^\top = A-ZAZ^\top-QAZ^\top-ZAQ^\top-QAQ^\top$. Let us define $\Delta=QAZ^\top+ZAQ^\top+QAQ^\top$.
Then, using the outer product form of $Q$, we get $\Delta =e_1\bigl(e_n^\top AZ^\top +(e_n^\top Ae_n)e_1^\top\bigr) +(ZAe_n)e_1^\top$. So we know that $\Delta$ is a sum of two outer products, each of inner rank at most one. By Lemma~\ref{lem:inner-rank-properties}, we have $\rho(\Delta)\leq 2$.

Since $\alpha_+(A)=\rho(A-ZAZ^\top)$, subadditivity
of inner rank gives $\alpha_c(A) \leq\alpha_+(A)+2$,
where we know $\rho(-\Delta)=\rho(\Delta)$ because $-I_n$ is invertible.

For the lower bound we rewrite the same identity as $A-ZAZ^\top=(A-Z_cAZ_c^\top)+\Delta$. Then, applying subadditivity again gives $\alpha_+(A) \leq\alpha_c(A)+\rho(\Delta) \leq\alpha_c(A)+2$.
Finally, rearranging and combining the two inequalities we conclude $\alpha_+(A)-2\leq\alpha_c(A)\leq\alpha_+(A)+2$, as desired.
\end{proof}

So, replacing the nilpotent shift by the cyclic shift
changes the displacement rank by at most two. However,
the corresponding displacement equation does not need to have a solution for every right-hand side. We next characterize the right-hand sides for which a solution exists.

\begin{lemma}
\label{lem:cyclic-displacement-solvability}
Let ${R}$ be a commutative ring with identity,
let $Z_c$ be the $n\times n$ lower cyclic shift matrix, and let $G\in{R}^{n\times n}$. Then the equation $A-Z_cAZ_c^\top=G$ has a solution if and only if the sum of the entries of $G$ along each cyclic diagonal is zero.
Equivalently, it has a solution if and only if
\[
\sum_{k=1}^{n}Z_c^{k-1}G(Z_c^\top)^{k-1}=0.
\]
\end{lemma}

\begin{proof}
Suppose first that $A-Z_cAZ_c^\top=G$.
Since $Z_c^n=I$, summing the shifted equations gives
\[
\sum_{k=1}^{n}Z_c^{k-1}G(Z_c^\top)^{k-1} = \sum_{k=1}^{n} \left(Z_c^{k-1}A(Z_c^\top)^{k-1} -Z_c^{k}A(Z_c^\top)^{k} \right) = 0.
\]
Each entry of this sum is the sum of the entries
of $G$ along the corresponding cyclic diagonal.
Thus, every cyclic diagonal sum is zero.

Conversely, suppose that each cyclic diagonal sum
is zero. We construct a solution entrywise.
Fix $d\in\{0,\dots,n-1\}$, and interpret column
indices cyclically modulo $n$, with values in
$\{1,\dots,n\}$.
Write $g_i=G_{i,i+d}$ for $1\leq i\leq n$.
By assumption, $\sum_{i=1}^n g_i=0$.

Define $s_1=0$ and
$s_i=\sum_{j=2}^i g_j$ for $2\leq i\leq n$,
and set $A_{i,i+d}=s_i$.
Then $s_i-s_{i-1}=g_i$ for $2\leq i\leq n$.
At the cyclic boundary, $s_1-s_n=-\sum_{j=2}^n g_j=g_1$.
Repeating this construction for every $d$ defines
all entries of $A$ and gives
$A_{ij}-A_{i-1,j-1}=G_{ij}$ with cyclic indexing.
Hence, $A-Z_cAZ_c^\top=G$, as desired.
\end{proof}

Now that we have characterized solvability, we construct all the solutions using the Toeplitz reconstruction from
Lemma \ref{lem:unique-sol-pir}. The boundary error defined below is the
correction computed in Lemma \ref{lem:comparison-nilpotent-cyclic-displacement}, evaluated at the
Toeplitz reconstruction $A_0$.

\begin{defn}
\label{def:boundary-error}
Let $R$ be a commutative ring, and let $x_\ell,y_\ell\in R^n$ for $\ell=1,\dots,\alpha$. Set $G=\sum_{\ell=1}^{\alpha}x_\ell y_\ell^\top$ and $A_0=\sum_{\ell=1}^{\alpha}L(x_\ell)U(y_\ell^\top)$, and let $\mathcal{L}_c(Y)=Y-Z_cYZ_c^\top$ denote the cyclic displacement operator. The \emph{boundary error matrix} is $E=G-\mathcal{L}_c(A_0)$. We also define, for $2\le i,j\le n$,
\[
a_j=\sum_{\ell=1}^{\alpha}\sum_{k=1}^{j-1}(x_\ell)_{n-k+1}(y_\ell)_{j-k},
\qquad
b_i=\sum_{\ell=1}^{\alpha}\sum_{k=1}^{i-1}(x_\ell)_{i-k}(y_\ell)_{n-k+1}.
\]
\end{defn}

\begin{lemma}
\label{lem:boundary-error}
With the notation of \cref{def:boundary-error}, all entries of $E$ outside its first row and first column are zero, and
\[
E=
\begin{bmatrix}
E_{11} & a_2 & a_3 & \cdots & a_n\\
b_2 & 0 & 0 & \cdots & 0\\
b_3 & 0 & 0 & \cdots & 0\\
\vdots & \vdots & \vdots & \ddots & \vdots\\
b_n & 0 & 0 & \cdots & 0
\end{bmatrix},
\qquad\text{where }
E_{11}=\sum_{\ell=1}^{\alpha}\sum_{m=1}^{n}(x_\ell)_m(y_\ell)_m .
\]
\end{lemma}

\begin{proof}
By \cref{lem:unique-sol-pir}, $A_0-ZA_0Z^\top=G$. Substituting $Z_c=Z+e_1e_n^\top$ into $\mathcal{L}_c(A_0)$ gives $\mathcal{L}_c(A_0)=(A_0-ZA_0Z^\top)-e_1e_n^\top A_0Z^\top-ZA_0e_ne_1^\top-e_1e_n^\top A_0e_ne_1^\top$, and therefore
\[
E=G-\mathcal{L}_c(A_0)=e_1(e_n^\top A_0Z^\top)+(ZA_0e_n)e_1^\top+(e_n^\top A_0e_n)e_1e_1^\top .
\]
Write $r^\top=e_n^\top A_0$ for the last row of $A_0$ and $c=A_0e_n$ for the last column.
For a row vector $r^\top$ we have $(r^\top Z^\top)_1=0$ and $(r^\top Z^\top)_j=r_{j-1}$ for $j\ge2$,
and for a column vector $c$ we have $(Zc)_1=0$ and $(Zc)_i=c_{i-1}$ for $i\ge2$.
Hence the three terms of $E$ are supported on the first row, the first column, and the
$(1,1)$ position, respectively, so $E$ vanishes outside its first row and first column, and that results in $E_{11}=(A_0)_{nn}$, $E_{1j}=(A_0)_{n,j-1}$ for $2\le j\le n$ and $E_{i1}=(A_0)_{i-1,n}$ for $2\le i\le n$.
To evaluate these, note that $L(x)_{ik}=x_{i-k+1}$ for $i\ge k$ and $U(y^\top)_{kj}=y_{j-k+1}$
for $j\ge k$, and both vanish otherwise, so we get $\bigl(L(x)U(y^\top)\bigr)_{ij}=\sum_{k=1}^{\min(i,j)}x_{i-k+1}\,y_{j-k+1}$.
Summing over $\ell$ gives us $(A_0)_{nn}= E_{11}$, $(A_0)_{n,j-1}=a_j$, and $(A_0)_{i-1,n}=b_i$ which are the entries displayed in the statement.
\end{proof}

\begin{lemma}
\label{lem:circulant-displacement}
With the notation of \cref{def:boundary-error}, the equation $A-Z_cAZ_c^\top=G$ has a solution in ${R}^{n\times n}$ if and only if $E_{11}=0$ and $b_i=-a_{n-i+2}$ for $2\leq i\leq n$. When these conditions hold, all solutions are given by $A=A_0+X+C$, where $C$ is an arbitrary circulant matrix and $X$ is the strictly upper-triangular Toeplitz matrix defined by
\[
X_{ij}=
\begin{cases}
a_{j-i+1}, & i<j,\\
0, & i\geq j.
\end{cases}
\]
\end{lemma}

\begin{proof}
Let $\mathcal{L}_c$ be the cyclic displacement operator of \cref{def:boundary-error}. Suppose $A_1$ and $A_2$ are two solutions of $\mathcal{L}_c(A)=G$, and define $\Delta=A_1-A_2$. Then, $\mathcal{L}_c(\Delta)=0$, which means $\Delta=Z_c\Delta Z_c^{\top}$.
Entrywise, this means $\Delta_{ij}=\Delta_{i-1,j-1}$ where indices are cyclic modulo $n$. Then, $\Delta$ is constant on the cyclic diagonals, so $\Delta$ is circulant. This means that our general solution is going to be of the form $A=A_p+C$, where $A_p$ is a particular solution and $C$ is an arbitrary circulant matrix.

We begin by constructing a particular solution. By the definition of $E$, we have $\mathcal{L}_c(A_0)=G-E$, and by \cref{lem:boundary-error}, $E$ vanishes outside its first row and first column, with first row $(E_{11},a_2,\dots,a_n)$ and first column $(E_{11},b_2,\dots,b_n)^\top$.

To prove necessity, we begin with the assumption that $A$ is a solution. Then we have that $\mathcal{L}_c(A-A_0)=G-(G-E)=E$, so by the necessity direction of Lemma \ref{lem:cyclic-displacement-solvability} applied with $E$ in place of $G$, every cyclic
diagonal sum of $E$ is zero. Since $E$ is supported on its first row and first column, we check which support positions lie on each cyclic diagonal $\{(i,j):j-i\equiv d \pmod n\}$. For $d=0$ the only support position is $(1,1)$ since $(1,j)$ with $j\equiv1$ forces $j=1$ and $(i,1)$ with $1\equiv i$ forces $i=1$. Hence $E_{11}=0$. Now, for $1\le d\le n-1$ the support positions are $(1,1+d)$ in the first row and $(n-d+1,1)$ in the first column, from $1-i\equiv d$, and no others. Hence $0=E_{1,1+d}+E_{n-d+1,1}=a_{d+1}+b_{n-d+1}$. Setting $i=n-d+1\in\{2,\dots,n\}$
gives $b_i=-a_{n-i+2}$.
Therefore $E$ has the form
\[
E=\begin{pmatrix}
0 & a_2 & \cdots & a_n\\
b_2 & 0 & \cdots & 0\\
\vdots & & & \vdots\\
b_n & 0 & \cdots & 0
\end{pmatrix}
\]
with $b_i=-a_{n-i+2}$. This proves necessity.

Conversely, suppose $E_{11}=0$ and $b_i=-a_{n-i+2}$ for
$2\le i\le n$. Then $E$ has the displayed form, and $X$ is defined as in the statement of the lemma. Now, we claim that $\mathcal{L}_c(X)=E$. Notice, $(\mathcal{L}_c(X))_{ij} = X_{ij}-X_{i-1,j-1}$ with cyclic indexing. If $i=1$ and $j>1$, then $X_{1j}-X_{n,j-1}=a_j-0=a_j$, we get the first row of $E$. Also, if $j=1$ and $i>1$, then $X_{i1}-X_{i-1,n}=0-a_{n-i+2}=b_i$, we get the first column of $E$. At $(i,j)=(1,1)$, we have
$X_{11}-X_{nn}=0=E_{11}$. For $i,j\geq 2$, the Toeplitz structure of $X$ gives $X_{ij}-X_{i-1,j-1}=0=E_{ij}$. Then, $X-Z_cXZ_c^{\top}=E$.

Hence, $\mathcal{L}_c(A_0+X) = \mathcal{L}_c(A_0) + \mathcal{L}_c(X) = (G-E)+E = G$, which shows that $A_0+X$ is a particular solution. Conversely, every circulant matrix $C$ satisfies $\mathcal{L}_c(C)=0$. So adding the arbitrary circulant matrix $C$ gives the general solution:
\[
A = \sum_{\ell=1}^{\alpha}L(x_\ell)U(y_\ell^{\top}) + X + C.
\]
\end{proof}

\subsection{Circulant Matrices and the Discrete Fourier Transform}

Lemma \ref{lem:circulant-displacement} shows that every solution of the cyclic displacement equation differs from the particular solution $A_0+X$ by an arbitrary circulant matrix $C$. To understand and invert such a
matrix, we diagonalize it. Over a ring containing a suitable root of unity, the discrete Fourier
transform does this. We first define the transform and show that it is invertible, then test the hypotheses on examples, then diagonalize circulant
matrices and deduce when they are invertible. These additional assumptions on the coefficient ring are not needed for the solvability and reconstruction results above.

\begin{defn}[Discrete Fourier Transform Matrix]
\label{def:discrete-fourier-transform-matrix}
Let $R$ be a commutative ring with identity, and let
$n\geq 1$ be an integer such that $n$ is a unit.
We write $n^{-1}$ for the inverse of $n$.
Let $\omega\in R$ be a principal $n$-th root of unity,
meaning that $\omega^n=1$ and
$\sum_{\ell=1}^{n}\omega^{t(\ell-1)}=0$
for $1\leq t\leq n-1$.
The $n\times n$ discrete Fourier transform matrix $F$
is defined by $F_{jk}=\omega^{(j-1)(k-1)}$
for $1\leq j,k\leq n$.
Since $\omega^n=1$, the element $\omega$ is a unit
with $\omega^{-1}=\omega^{n-1}$.
\end{defn}

The next lemma shows that, under these assumptions, $F$ is invertible with an explicit inverse, which is what allows us to use $F$ as a change of basis.

\begin{lemma}
\label{lem:dft-inverse}
Under the assumptions of Definition~\ref{def:discrete-fourier-transform-matrix},
the matrix $F$ is invertible, where we have
$(F^{-1})_{jk}=n^{-1}\omega^{-(j-1)(k-1)}$
for $1\leq j,k\leq n$.
\end{lemma}

\begin{proof}
Let $V\in R^{n\times n}$ be defined by
$V_{jk}=\omega^{-(j-1)(k-1)}$.
Then
$(FV)_{ij}=\sum_{k=1}^{n}\omega^{(i-j)(k-1)}$.
If $i=j$, this sum is $n\cdot 1_R$.
If $i\neq j$, the difference $i-j$ is congruent
modulo $n$ to some $t\in\{1,\dots,n-1\}$.
Since $\omega^n=1$, the principal-root condition gives
$(FV)_{ij}=\sum_{k=1}^{n}\omega^{t(k-1)}=0$.
Thus, $FV=nI_n$.
The same computation gives
$(VF)_{ij}=\sum_{k=1}^{n}\omega^{(j-i)(k-1)}$,
so $VF=nI_n$ as well.
Since $n$ is a unit,
$F(n^{-1}V)=(n^{-1}V)F=I_n$.
Hence $F^{-1}=n^{-1}V$, as desired.
\end{proof}

Before using $F$, we check how restrictive the hypotheses of Definition \ref{def:discrete-fourier-transform-matrix} are.
First, over $\mathbb{Z}$ with $n=2$,
the element $\omega=-1$ is a principal second root
of unity, but the corresponding Fourier matrix has
determinant $-2$ and is not invertible over $\mathbb{Z}$.
Thus, the principal-root condition alone does not
guarantee invertibility. This is why Definition \ref{def:discrete-fourier-transform-matrix} also requires $n$ to be a unit.

Next, over a finite field $\mathbb{F}_q$ of characteristic
$p$, these assumptions hold for some $\omega$ exactly
when $n\mid q-1$.
Indeed, $\mathbb{F}_q^\times$ is cyclic of order $q-1$,
so this divisibility condition gives an element
$\omega$ of order $n$ and implies $p\nmid n$.
For $1\leq t<n$, we have $1-\omega^t\neq 0$, and
$(1-\omega^t)\sum_{\ell=1}^{n}\omega^{t(\ell-1)}
=1-\omega^{tn}=0$.
Since $\mathbb{F}_q$ is a field, each sum is zero.
Conversely, if the assumptions hold and $\omega$
had order $d<n$, then
$\sum_{\ell=1}^{n}\omega^{d(\ell-1)}=n\cdot 1$
would be both zero and a unit, a contradiction.
Thus, $\omega$ has order $n$, so $n\mid q-1$.

Finally, the assumptions also hold over rings with zero divisors.
For example, take $R=\mathbb{Z}/25\mathbb{Z}$,
$n=4$, and $\omega=7$.
Then, $4$ is a unit and $\omega^4=1$.
For $t=1,2,3$, the elements $1-\omega^t$ are units.
The same geometric-sum identity therefore shows that
$\omega$ is a principal fourth root of unity.

With these examples in hand, we return to the circulant matrices of Lemma \ref{lem:circulant-displacement}. Since a circulant
matrix is a polynomial in the cyclic shift $Z_c$, it suffices to diagonalize $Z_c$ by $F$.

\begin{theorem}
\label{thm:DFT-diagonalization}
Under the assumptions of Definition \ref{def:discrete-fourier-transform-matrix}, every circulant
matrix $C$ over $R$ is diagonalized by the
discrete Fourier transform matrix $F$. In other words, if $C=\sum_{k=1}^{n}c_kZ_c^{k-1}$, then
$FCF^{-1}=D_C$, where $D_C$ is diagonal and its entries are defined by $(D_C)_{jj}=\sum_{k=1}^{n}c_k\omega^{(j-1)(k-1)}$ for $1\leq j\leq n$.
\end{theorem}

\begin{proof}
Since $C$ is a polynomial in $Z_c$, we first show that $F$ diagonalizes $Z_c$. For each $j=1,\dots,n$, define the column vector $v^{(j)}$ by $v^{(j)}_k=\omega^{-(j-1)(k-1)}$ for $1\leq k\leq n$.

For $2\leq k\leq n$,
$(Z_cv^{(j)})_k=v_{k-1}^{(j)}
=\omega^{j-1}v_k^{(j)}$.
For the first component,
$(Z_cv^{(j)})_1=v_n^{(j)}
=\omega^{-(j-1)(n-1)}
=\omega^{j-1}v_1^{(j)}$.

Let $V=[v^{(1)}\;\cdots\;v^{(n)}]$. By Lemma~\ref{lem:dft-inverse}, we know $V=nF^{-1}$. So, $V$ is invertible, and the identities give $Z_cV=VD_Z$, where $D_Z=\operatorname{diag}(1,\omega,\dots,\omega^{n-1})$. Then, $V^{-1}Z_cV=D_Z$. Since $V^{-1}=n^{-1}F$, this is equivalent to $FZ_cF^{-1}=D_Z$. Then, we can apply the similarity transformation to the full circulant matrix. We get
\[
FCF^{-1}=\sum_{k=1}^{n}c_kD_Z^{k-1}.
\]

Since sums and powers of diagonal matrices stay diagonal, this matrix is diagonal. Its $j$-th diagonal entry is $(D_C)_{jj}=\sum_{k=1}^{n}c_k\omega^{(j-1)(k-1)}$, which is the $j$-th component of the discrete Fourier transform of $(c_1,\dots,c_{n})^\top$ under our convention.
\end{proof}

Theorem \ref{thm:DFT-diagonalization} reduces questions about a circulant matrix to questions about its eigenvalues
$\lambda_1,\dots,\lambda_n\in R$. In particular, it gives a criterion for invertibility and a formula for the inverse.

\begin{cor}
\label{cor:circulant-inverse}
Let ${R}$ be a commutative ring with identity,
assume that $n\cdot 1_{{R}}$ is a unit,
and let $\omega$ be a principal $n$-th root of unity.
Let $F$ be the corresponding discrete Fourier transform matrix, and let
$C=\sum_{k=1}^{n}c_kZ_c^{k-1}$ be a circulant matrix.
Define $\lambda_j=\sum_{k=1}^{n}c_k
\omega^{(j-1)(k-1)}$ for $1\leq j\leq n$. Then $C$ is invertible if and only if every $\lambda_j$ is a unit in ${R}$. In that case, $C^{-1} = F^{-1} \operatorname{diag}(\lambda_1^{-1},\dots,\lambda_{n}^{-1}) F$, and $C^{-1}$ is also circulant.
\end{cor}

\begin{proof}
By \cref{thm:DFT-diagonalization}, $FCF^{-1}=D_C$, where
$D_C=\operatorname{diag}(\lambda_1,\dots,\lambda_{n})$.
Since $F$ is invertible, $C$ is invertible if and
only if $D_C$ is invertible.

If every $\lambda_j$ is a unit, then
$D_C^{-1}=\operatorname{diag}
(\lambda_1^{-1},\dots,\lambda_{n}^{-1})$.
Conversely, if $D_C$ has an inverse $B$, the
diagonal entries of $D_CB=I$ give
$\lambda_jB_{jj}=1$, so every $\lambda_j$ is a unit.
Taking inverses in $C=F^{-1}D_CF$ gives the
stated formula for $C^{-1}$.

Finally, since $C$ is a polynomial in $Z_c$,
we have $CZ_c=Z_cC$.
Multiplying on the left and right by $C^{-1}$
gives $Z_cC^{-1}=C^{-1}Z_c$.
Since $Z_c^{-1}=Z_c^\top$, this implies $C^{-1}=Z_cC^{-1}Z_c^\top$. Entrywise, this says that $C^{-1}$ is constant along each cyclic diagonal. Hence $C^{-1}$ is circulant.
\end{proof}

Corollary \ref{cor:circulant-inverse} relies on the hypotheses of Definition \ref{def:discrete-fourier-transform-matrix}. Over a field, invertibility can also be
decided without them, as follows.
Let $\mathbb K$ be a field and let
$C=\sum_{k=1}^{n}c_kZ_c^{k-1}$.
Define $c(x)=\sum_{k=1}^{n}c_kx^{k-1}$.
Then $C$ is invertible if and only if
$\gcd(c(x),x^n-1)=1$ in $\mathbb K[x]$.
If the greatest common divisor is $1$,
Bezout's identity gives polynomials $u,v\in \mathbb K[x]$
such that $u(x)c(x)+v(x)(x^n-1)=1$.
Evaluating at $Z_c$ gives $u(Z_c)C=Cu(Z_c)=I_n$.
Conversely, if $C$ is invertible, its inverse commutes
with $Z_c$ and is therefore circulant, as shown in
the previous proof. Write $C^{-1}=u(Z_c)$.
The remainder of $u(x)c(x)-1$ upon division by
$x^n-1$ has degree less than $n$ and vanishes at $Z_c$.
Since $I_n,Z_c,\dots,Z_c^{n-1}$ are linearly independent, this remainder
is zero. Hence $u(x)c(x)\equiv 1\pmod{x^n-1}$,
which proves the criterion.
This argument applies in every characteristic,
including when the characteristic divides $n$.

Finally, we return to the cyclic displacement equation of Section \ref{ss:4.2}. Under the hypotheses of
Definition \ref{def:discrete-fourier-transform-matrix}, the circulant part $C$ in Lemma \ref{lem:circulant-displacement} equals $F^{-1}D_CF$ by Theorem \ref{thm:DFT-diagonalization}, so every solution of $A-Z_cAZ_c^\top=G$ has the form $A_0+X+F^{-1}D_CF$, where $D_C$ is diagonal.

\section{Levinson-Type Recursions for Toeplitz Systems}
\label{sec:levinson}

This section turns from displacement representations to
recursive computations for Toeplitz matrices. The recursions use the repeated leading and trailing submatrices of a Toeplitz matrix to extend solutions one step at a time. Over a commutative ring, we must also show that every quantity inverted in these updates is invertible. We first establish recursions for symmetric block Toeplitz matrices, then treat nonsymmetric scalar Toeplitz systems, and finally bound the computational cost. Along the way, Lemma \ref{lem:denominator} expresses the recursion denominator in terms of leading principal determinants.


\subsection{Symmetric Block Toeplitz Matrices}

Let $R$ be a commutative ring and fix $p\ge 1$. For blocks $r_k\in R^{p\times p}$,
$-(s-1)\le k\le s-1$, let $\Gamma_m=[r_{i-j}]_{i,j=1}^{m}\in R^{mp\times mp}$ for $1\le m\le s$. This is the $m$-th leading principal block submatrix of $\Gamma_s$. The matrix $\Gamma_s$ is
\emph{block Toeplitz}, since its blocks are constant along each block diagonal, and it is
\emph{symmetric} if $\Gamma_s^\top=\Gamma_s$, equivalently $r_{-k}=r_k^\top$ for every $k$. For a block matrix
$W=\begin{pmatrix}W_{11}&W_{12}\\ W_{21}&W_{22}\end{pmatrix}$ with $W_{11}$ invertible, the Schur
complement of $W_{11}$ in $W$ is $W_{22}-W_{21}W_{11}^{-1}W_{12}$, and
$\det W=\det W_{11}\,\det\bigl(W_{22}-W_{21}W_{11}^{-1}W_{12}\bigr)$. Symmetrically, if $W_{22}$ is invertible, then we know the Schur complement of $W_{22}$ in $W$ is $W_{11}-W_{12}W_{22}^{-1}W_{21}$ and
$\det W=\det W_{22}\,\det\bigl(W_{11}-W_{12}W_{22}^{-1}W_{21}\bigr)$.

The block recursions require inverses of the residual
matrices $M_m$ and $N_m$. We first show that these matrices are invertible whenever the leading principal block submatrices are invertible.

\begin{lemma}
\label{lem:levinson-residual-invertibility}
Let $R$ be a commutative ring and let $\Gamma_s=[r_{i-j}]_{i,j=1}^s$ be a symmetric block Toeplitz
matrix with $p\times p$ blocks. Assume $\Gamma_m$ is invertible for $1\le m\le s$. Suppose that for
some $m$ with $1\le m\le s$ there are matrices $A_m,B_m\in R^{mp\times p}$ and
$M_m,N_m\in R^{p\times p}$ with
\[
\Gamma_mA_m=\begin{pmatrix}0\\ \vdots\\ 0\\ N_m\end{pmatrix},\qquad
\Gamma_mB_m=\begin{pmatrix}M_m\\ 0\\ \vdots\\ 0\end{pmatrix},
\]
where the last block of $A_m$ and the first block of $B_m$ equal $I_p$. Then $N_m$ and $M_m$ are invertible.

\end{lemma}

\begin{proof}
For $m=1$, the normalization gives $A_1=B_1=I_p$,
so $N_1=M_1=r_0=\Gamma_1$, which is invertible by
assumption.

For $m\ge2$, consider the equation $\Gamma_mA_m=\left[\begin{smallmatrix}0\\ \vdots\\ 0\\ N_m\end{smallmatrix}\right]$.
Then, we will write $\Gamma_m=\begin{pmatrix}\Gamma_{m-1}&K_m\\ K_m^\top&r_0\end{pmatrix}$, where $K_m$ is the block
column with entries $r_{1-m},\dots,r_{-1}$. Note that the lower-left block is $K_m^\top$ because
$r_{-k}^\top=r_k$. By the Toeplitz structure the lower-right $(m-1)p\times(m-1)p$ block of $\Gamma_m$
is also $\Gamma_{m-1}$, so also
$\Gamma_m=\begin{pmatrix}r_0&H_m\\ H_m^\top&\Gamma_{m-1}\end{pmatrix}$, where $H_m$ is the block row
with entries $r_{-1},\dots,r_{1-m}$. We also write $A_m=\begin{pmatrix}\widehat A_m\\ I_p\end{pmatrix}$.
Multiplying gives
$\Gamma_{m-1}\widehat A_m+K_m=0$ and
$K_m^\top\widehat A_m+r_0=N_m$.
Since $\Gamma_{m-1}$ is invertible, we will get
$\widehat A_m=-\Gamma_{m-1}^{-1}K_m$.
Substituting into the second equation gives us $N_m=r_0-K_m^\top \Gamma_{m-1}^{-1}K_m$. Thus, $N_m$ is the Schur complement of $\Gamma_{m-1}$ in $\Gamma_m$. The Schur complement determinant identity gives $\det(\Gamma_m)=\det(\Gamma_{m-1})\det(N_m)$. Since $\det(\Gamma_m)$ and $\det(\Gamma_{m-1})$ are units, $\det(N_m)$ is also a unit. Hence $N_m$ is invertible.

We can prove the invertibility of $M_m$ using a similar Schur complement argument. Write
$B_m=\begin{pmatrix}I_p\\ \widehat B_m\end{pmatrix}$.
The equation $\Gamma_mB_m=\left[\begin{smallmatrix}M_m\\ 0\\ \vdots\\ 0\end{smallmatrix}\right]$ gives
$H_m^\top+\Gamma_{m-1}\widehat B_m=0$ and $r_0+H_m\widehat B_m=M_m$.
Since $\Gamma_{m-1}$ is invertible, $\widehat B_m=-\Gamma_{m-1}^{-1}H_m^\top$, and therefore
$M_m=r_0-H_m\Gamma_{m-1}^{-1}H_m^\top$. Thus $M_m$ is the Schur complement of the lower-right
block $\Gamma_{m-1}$ in $\Gamma_m$, and the second determinant identity above gives
$\det(\Gamma_m)=\det(\Gamma_{m-1})\det(M_m)$. Again, both $\det(\Gamma_m)$ and
$\det(\Gamma_{m-1})$ are units, so $\det(M_m)$ is a unit. Thus $M_m$ is invertible, as desired.
\end{proof}

\begin{theorem}
\label{thm:generalized-levinson-algorithm}
Let $R$ be a commutative ring with identity, and let
$\Gamma_s=[r_{i-j}]_{i,j=1}^s$ be a symmetric block-Toeplitz
matrix, where the blocks are $p\times p$ matrices
over $R$ satisfying $r_{-i}=r_i^\top$.
Assume that every leading principal block submatrix
$\Gamma_m$ is invertible for $m=1,\dots,s$.

Then matrices $A_m,B_m\in R^{mp\times p}$ and
invertible matrices $M_m,N_m\in R^{p\times p}$
satisfying
\[
\Gamma_mA_m=
\begin{bmatrix}
0\\ \vdots\\ 0\\ N_m
\end{bmatrix},
\qquad
\Gamma_mB_m=
\begin{bmatrix}
M_m\\ 0\\ \vdots\\ 0
\end{bmatrix},
\]
with the last block of $A_m$ and the first block
of $B_m$ equal to $I_p$, can be computed using
the following generalized Levinson recursions:
\begin{align*}
A_{m+1}
&=
\begin{bmatrix}0\\ A_m\end{bmatrix}
-
\begin{bmatrix}B_m\\ 0\end{bmatrix}
M_m^{-1}\Delta_m^\top,\\
B_{m+1}
&=
\begin{bmatrix}B_m\\ 0\end{bmatrix}
-
\begin{bmatrix}0\\ A_m\end{bmatrix}
N_m^{-1}\Delta_m,\\
N_{m+1}
&=N_m-\Delta_mM_m^{-1}\Delta_m^\top,\\
M_{m+1}
&=M_m-\Delta_m^\top N_m^{-1}\Delta_m,
\end{align*}
for $1\leq m<s$, with
$A_1=B_1=I_p$, $M_1=N_1=r_0$, and $\Delta_m=
\begin{bmatrix}
r_m&r_{m-1}&\cdots&r_1
\end{bmatrix}B_m$.
\end{theorem}

\begin{proof}
We will use induction. The base case holds since
$\Gamma_1=r_0$, $A_1=B_1=I_p$, and $M_1=N_1=r_0$.

Assume that at step $m$,
\[
\Gamma_mA_m=
\begin{bmatrix}
0\\ \vdots\\ 0\\ N_m
\end{bmatrix},
\qquad
\Gamma_mB_m=
\begin{bmatrix}
M_m\\ 0\\ \vdots\\ 0
\end{bmatrix},
\]
and that the last block of $A_m$ and the first
block of $B_m$ are $I_p$.
By Lemma~\ref{lem:levinson-residual-invertibility},
$M_m$ and $N_m$ are invertible, so their inverses
are well-defined.

Let
\[
a=\begin{bmatrix}0\\ A_m\end{bmatrix},
\qquad
b=\begin{bmatrix}B_m\\ 0\end{bmatrix},
\qquad
Q_m=
\begin{bmatrix}
r_{-1}&r_{-2}&\cdots&r_{-m}
\end{bmatrix}A_m.
\]
Using the Toeplitz structure of $\Gamma_{m+1}$ gives
\[
\Gamma_{m+1}a=
\begin{bmatrix}
Q_m\\ 0\\ \vdots\\ 0\\ N_m
\end{bmatrix},
\qquad
\Gamma_{m+1}b=
\begin{bmatrix}
M_m\\ 0\\ \vdots\\ 0\\ \Delta_m
\end{bmatrix}.
\]
Since the last block of $A_m$ and the first block
of $B_m$ are $I_p$, we have
$a^\top \Gamma_{m+1}b=\Delta_m$ and
$(\Gamma_{m+1}a)^\top b=Q_m^\top$.
The symmetry of $\Gamma_{m+1}$ shows that these are
equal. Hence $Q_m=\Delta_m^\top$.

To find the recursion for $A_{m+1}$, we take a
linear combination of the shifted terms to
eliminate $\Delta_m^\top$:
\[
\Gamma_{m+1}
\left(
\begin{bmatrix}0\\ A_m\end{bmatrix}
-
\begin{bmatrix}B_m\\ 0\end{bmatrix}
M_m^{-1}\Delta_m^\top
\right)
=
\begin{bmatrix}
0\\ \vdots\\ 0\\
N_m-\Delta_mM_m^{-1}\Delta_m^\top
\end{bmatrix}.
\]
This gives the stated updates for $A_{m+1}$ and
$N_{m+1}$.

Similarly, to find $B_{m+1}$, we eliminate
$\Delta_m$:
\[
\Gamma_{m+1}
\left(
\begin{bmatrix}B_m\\ 0\end{bmatrix}
-
\begin{bmatrix}0\\ A_m\end{bmatrix}
N_m^{-1}\Delta_m
\right)
=
\begin{bmatrix}
M_m-\Delta_m^\top N_m^{-1}\Delta_m\\
0\\ \vdots\\ 0
\end{bmatrix}.
\]
This gives the stated updates for $B_{m+1}$ and
$M_{m+1}$.

The last block of $A_{m+1}$ and the first block
of $B_{m+1}$ remain $I_p$.
Since $\Gamma_{m+1}$ is invertible by assumption,
Lemma~\ref{lem:levinson-residual-invertibility}
also shows that $N_{m+1}$ and $M_{m+1}$ are
invertible. This completes the induction. Thus, the recursions are well-defined at every
step and use only inverses of invertible matrices
over $R$.
\end{proof}

\subsection{Nonsymmetric Toeplitz Systems}
We next consider scalar Toeplitz systems without a symmetry assumption. Over a field the
recursion only needs its denominators to be nonzero; over a commutative ring we must show
that they are units. Let $R$ be a commutative ring, let $T_n=[t_{i-j}]_{i,j=1}^n\in R^{n\times n}$ be a Toeplitz
matrix (so $t_j\in R$ for $-(n-1)\le j\le n-1$), and let $b\in R^n$. Write $T_k=[t_{i-j}]_{i,j=1}^k$ for the leading principal
submatrices of $T_n$, and assume that every $T_k$, $1\le k\le n$, is invertible. For
$1\le k\le n$ we define $u_k=T_k^{-1}e_1$, $v_k=T_k^{-1}e_k$, and $x_k=T_k^{-1}b_{1:k}$,
where $b_{1:k}$ consists of the first $k$ entries of $b$; in particular $x_n=T_n^{-1}b$.
For $1\le k<n$ we further define $w_k^\top=(t_k,\dots,t_1)$, $\gamma_k=w_k^\top u_k$,
$\sigma_k=(t_{-1},\dots,t_{-k})\,v_k$, $d_k=1-\gamma_k\sigma_k$, as well as $U_k=\begin{pmatrix}u_k\\0\end{pmatrix}$ and $V_k=\begin{pmatrix}0\\v_k\end{pmatrix}\in R^{k+1}$.
Write $T_{k+1}=\begin{pmatrix}T_k&c_k\\ w_k^\top&t_0\end{pmatrix}$ with
$c_k=(t_{-k},\dots,t_{-1})^\top$. Since $T_{k+1}$ is Toeplitz, its lower-right $k\times k$
block is also $T_k$, and its first row is $(t_0,t_{-1},\dots,t_{-k})$. Hence $T_{k+1}U_k=e_1+\gamma_ke_{k+1}$ and $T_{k+1}V_k=\sigma_ke_1+e_{k+1}$.

\begin{lemma}
\label{lem:denominator}
Let $R$, $T_n$ and $b$ be as above, with every leading principal submatrix $T_k$, $1\le k\le n$,
invertible. With $\det T_0=1$, for $1\le k<n$ we have
$d_k=\det(T_{k+1})\det(T_{k-1})\bigl(\det T_k\bigr)^{-2}$. In particular $d_k\in R^\times$.
\end{lemma}

\begin{proof}
Let $P_k=[\,U_k\ V_k\,]\in R^{(k+1)\times2}$, $E_k=[\,e_1\ e_{k+1}\,]$, and
$B_k=\begin{pmatrix}1&\sigma_k\\ \gamma_k&1\end{pmatrix}$, so $\det B_k=d_k$.
The identities $T_{k+1}U_k = e_1+\gamma_k e_{k+1}$ and
$T_{k+1}V_k = \sigma_k e_1 + e_{k+1}$ stated at the start of this subsection say that
the columns of $T_{k+1}P_k$ are the columns of $E_kB_k$, that is,
$T_{k+1}P_k = E_kB_k$. Hence $P_k = WB_k$ where $W = T_{k+1}^{-1}E_k$. Let $S\in R^{2\times(k+1)}$ have the rows $e_1^\top,e_{k+1}^\top$.
The last entry of $U_k$ and the first entry of $V_k$ are $0$, so then $SP_k=\begin{pmatrix}(u_k)_1&0\\0&(v_k)_k\end{pmatrix}$ and also $SW=N=\bigl(T_{k+1}^{-1}\bigr)_{\{1,k+1\},\{1,k+1\}}$.
Taking determinants of $SP_k=N B_k$ gives
\begin{equation}\label{eq:detN}
(u_k)_1(v_k)_k=\det(N)\,d_k .
\end{equation}
By the cofactor formula, $(u_k)_1=(T_k^{-1})_{11}$ and $(v_k)_k=(T_k^{-1})_{kk}$ both equal
$\det(T_{k-1})\det(T_k)^{-1}$, since deleting the first row and column, or the last row
and column, of $T_k$ leaves $T_{k-1}$. Jacobi's complementary minor identity says for an invertible matrix $M$ and index sets $I,J$ of equal size,
$\det\bigl((M^{-1})_{I,J}\bigr)=(-1)^{\sum I+\sum J}\det(M_{J^c,I^c})\,/\det M$. To justify this identity over a commutative ring, let
$r=|I|=|J|\geq 1$ and use the polynomial identity
$\det((\operatorname{adj}M)_{I,J})
=(-1)^{\sum I+\sum J}(\det M)^{r-1}
\det(M_{J^c,I^c})$.
Since $\det M$ is a unit and
$M^{-1}=(\det M)^{-1}\operatorname{adj}M$,
division by $(\det M)^r$ gives the stated identity.
Here $I=J=\{1,k+1\}$, with sign $+1$ and $J^c=\{2,\dots,k\}$. For $k=1$ the set is empty with
minor $1=\det T_0$. Hence, $\det N=\det\bigl((T_{k+1})_{\{2,\dots,k\},\{2,\dots,k\}}\bigr)\det(T_{k+1})^{-1}
=\det(T_{k-1})\det(T_{k+1})^{-1}$, since the middle block of $T_{k+1}$ is $T_{k-1}$.
All of these determinants are units, so $\det N\in R^\times$, and solving \eqref{eq:detN} for $d_k$ gives the claimed formula and $d_k\in R^\times$.
\end{proof}

\begin{theorem}
\label{thm:nonsym}
Let $R$, $T_n$ and $b$ be as above, with every $T_k$, $1\le k\le n$, invertible. Then, for $1\le k<n$, $d_k$ is a unit and $u_{k+1}=d_k^{-1}\bigl(U_k-\gamma_kV_k\bigr)$ and $ v_{k+1}=d_k^{-1}\bigl(V_k-\sigma_kU_k\bigr)$.
Moreover, with $\eta_k=b_{k+1}-w_k^\top x_k$, we have
$x_{k+1}=\begin{pmatrix}x_k\\0\end{pmatrix}+\eta_kv_{k+1}$.
Consequently we have $T_{k+1}u_{k+1}=e_1$, $T_{k+1}v_{k+1}=e_{k+1}$ and also
$T_{k+1}x_{k+1}=b_{1:k+1}$.
\end{theorem}

\begin{proof}
By Lemma~\ref{lem:denominator}, $d_k\in R^\times$. By the identities $T_{k+1}U_k = e_1+\gamma_k e_{k+1}$ and
$T_{k+1}V_k = \sigma_k e_1 + e_{k+1}$ stated previously, we get $T_{k+1}(U_k-\gamma_kV_k)=e_1+\gamma_ke_{k+1}-\gamma_k(\sigma_ke_1+e_{k+1})=d_ke_1$ and also $T_{k+1}(V_k-\sigma_kU_k)=\sigma_ke_1+e_{k+1}-\sigma_k(e_1+\gamma_ke_{k+1})=d_ke_{k+1}$,
using $\sigma_k\gamma_k=\gamma_k\sigma_k$. Multiplying by $d_k^{-1}$ gives $T_{k+1}u_{k+1}=e_1$ and $T_{k+1}v_{k+1}=e_{k+1}$. Here $u_{k+1}$ on the left is the vector $d_k^{-1}(U_k - \gamma_kV_k)$ from the statement. Now, since $T_{k+1}$ is invertible,
these are exactly $u_{k+1}=T_{k+1}^{-1}e_1$ and $v_{k+1}=T_{k+1}^{-1}e_{k+1}$. For the
last claim, let $X_k=\begin{pmatrix}x_k\\0\end{pmatrix}$. Then
we get that $T_{k+1}X_k=\begin{pmatrix}b_{1:k}\\ w_k^\top x_k\end{pmatrix}$, so $T_{k+1}(X_k+\eta_kv_{k+1})=\begin{pmatrix}b_{1:k}\\ w_k^\top x_k+\eta_k\end{pmatrix}=b_{1:k+1}.$
\end{proof}

The assumption that every leading principal submatrix is invertible cannot be replaced by invertibility of the full matrix in the stated recursion. For example, over $\mathbb{Z}/6\mathbb{Z}$, consider the Toeplitz matrix $T=\left[\begin{smallmatrix}2&1\\3&2\end{smallmatrix}\right]$. Its determinant is $2\cdot 2-1\cdot 3=1$, so $T$ is invertible. However, its first leading principal submatrix is $T_1=[2]$, and $2$ is not a unit modulo $6$. Thus, not every leading principal submatrix of $T$ is invertible, and the initialization $u_1=v_1=2^{-1}$ is not defined. This example shows why the leading-principal invertibility assumption is needed for these recursions.

Theorem \ref{thm:generalized-levinson-algorithm} constructs the residuals and the vectors $A_m,B_m$ for symmetric block Toeplitz
matrices but does not solve a linear system, and we do not give an operation count for it. Lemma \ref{lem:denominator} and Theorem \ref{thm:nonsym} treat nonsymmetric Toeplitz systems only in the scalar case. We do not treat the nonsymmetric block case.

\subsection{Operation Count}

We conclude by counting the ring operations required to solve the system.

\begin{cor}
\label{cor:levinson-computational-cost}
Under the assumptions of Theorem~\ref{thm:nonsym}, the recursions
find the solution $x_n=T_n^{-1}b$ using $O(n^2)$ additions and multiplications in $R$, at most $n$ inversions of units in $R$, and $O(n)$ ring elements of working storage beyond the input.
These bounds count ring operations rather than bit operations.
\end{cor}

\begin{proof}
We begin with $u_1=v_1=t_0^{-1}$ and $x_1=t_0^{-1}b_1$, which require one inversion of the unit $t_0$ and a constant number of other ring operations.

At step $k$, computing $\gamma_k$, $\sigma_k$, and $\eta_k$ requires three dot products of length $k$, so this takes $O(k)$ additions and multiplications. Computing $d_k$ takes a constant number of operations, and its inverse requires one
inversion of a unit.

The updates for $u_{k+1}$, $v_{k+1}$, and $x_{k+1}$ each involve vectors of length $k+1$, so they also require $O(k)$ additions and multiplications. Summing over $k=1,\dots,n-1$, the total is
\[
O\left(1+\sum_{k=1}^{n-1}k\right)=O(n^2).
\]
There is one initial inversion and at most one inversion at each subsequent step, giving at most $n$ inversions of units.

Finally, each step requires only the current vectors $u_k$, $v_k$, and $x_k$, together with space for their updates. Previous vectors can be discarded after each step, so the working storage is $O(n)$ ring elements.
\end{proof}

\bibliographystyle{amsalpha}
\bibliography{ref.bib}

\end{document}